\documentclass[letterpaper, twocolumn, 10pt, conference]{ieeeconf}
\usepackage[utf8]{inputenc}

\IEEEoverridecommandlockouts

\usepackage{amsmath}
\usepackage{amssymb}
\usepackage{amsthm}

\usepackage{bm}
\usepackage{xfrac}
\usepackage[ruled, vlined, linesnumbered]{algorithm2e}
\usepackage{tikz}
\usepackage{subcaption}
\usepackage{blindtext}

\usepackage{color}

\newcommand{\cK}{\mathcal{K}}

\newcommand{\Gu}{\mathbf{G}_{\bm{u}}}
\newcommand{\Gw}{\mathbf{G}_{\bm{w}}}
\renewcommand{\t}{^\mathrm{T}}

\newtheorem{proposition}{Proposition}

\title{ Constrained Covariance Steering Based Tube-MPPI
}

\author{Isin M. Balci \and Efstathios Bakolas \and Bogdan Vlahov \and Evangelos A. Theodorou   \thanks{This research has been supported in part by NSF awards CMMI-1937957, CMMI-1936079, and ECCS-1924790.
I. M. Balci (PhD student) and E. Bakolas (Associate Professor) are with the Department of Aerospace Engineering and Engineering Mechanics, The University of Texas at Austin, Austin, Texas 78712-1221, USA, Email: isinmertbalci@utexas.edu, bakolas@austin.utexas.edu. B. Vlahov (PhD student) and E. Theodorou (Associate Professor) are with the Department of Aerospace Engineering, Georgia Tech, Atlanta, Georgia, USA, Email: evangelos.theodorou@gatech.edu}}

\begin{document}

\maketitle

\begin{abstract}
    In this paper, we present a new trajectory optimization algorithm for stochastic linear systems which combines Model Predictive Path Integral (MPPI) control with Constrained Covariance Steering (CSS) to achieve high performance with safety guarantees (robustness). Although MPPI can be used to solve complex nonlinear trajectory optimization problems, it may not always handle constraints effectively and its performance may degrade in the presence of unmodeled disturbances. By contrast, CCS can handle probabilistic state and / or input constraints (e.g., chance constraints) by controlling uncertainty which implies that CCS can provide robustness against stochastic disturbances. CCS, however, suffers from scalability issues and cannot handle complex cost functions in general. We argue that the combination of the two methods yields a class of trajectory optimization algorithms that can achieve high performance while ensuring safety with high probability. The efficacy of our algorithm is demonstrated in an obstacle avoidance problem and a path generation problem with a circular track.
\end{abstract}

\section{Introduction}\label{sec:introduction}
Many real-world tasks for autonomous systems can be cast as finite horizon \textit{stochastic} trajectory optimization problems in the presence of model uncertainties, and random exogenous inputs from the environment. The main goal of these problems is to find control policies that minimize the expected value of a given cost function while satisfying state and input constraints with a given confidence level. 


In this work, we present a novel algorithm for constrained stochastic trajectory optimization problems subject to safety constraints. 
Our proposed algorithm combines Constrained Covariance Steering (CCS) theory for discrete-time stochastic linear systems with Model Predictive Path Integral (MPPI) to achieve robustness to uncertainties and variations of the different parameters of the proposed controllers
as well as improved performance, and scalability.

\noindent\textbf{\text{Literature Review:}}
Optimization-based methods treat the stochastic trajectory optimization problem as a nonlinear program (NLP) which can be solved by specialized NLP solvers. 
However, these NLP based approaches rely on a good initial guess to achieve high performance and may suffer from the lack of convergence guarantees \cite{p:gros2020nonlinearMPC}. 
Successive convexification-based methods provide convergence guarantees, but they may still suffer from scalability issues, if the underlying system dynamics are stochastic \cite{p:malyuta2021convexGeneration, mao2016successive, p:dandrea2012generationscp}.

Dynamic programming-based algorithms have been proposed for unconstrained stochastic trajectory optimization problems to alleviate the scalability issue in \cite{p:theodorou2010sddp, p:todorov2005iterativeLQG} and for constrained problems in \cite{p:aoyama2020constrainedddp, p:xie2017constrainedddp}. These methods, however, lack safety guarantees and their applicability is limited to smooth objective functions. On the other hand, sampling-based stochastic optimization algorithms deal with non-smooth objective functions but they often cannot handle unmodeled disturbances and model mismatch \cite{p:williams2017informationMPPI, p:williams2018tubeMPPI}. 

The two methods that are most closely related to our approach are discussed in  \cite{p:okamoto2019pathplanningcovariance}, \cite{p:yin2021improvingMPPIwithCS}. In \cite{p:okamoto2019pathplanningcovariance}, the authors use Covariance Steering for path planning for linear systems with chance constraints for obstacle avoidance.
However, the approach in \cite{p:okamoto2019pathplanningcovariance} does not scale well due to the fact that the number of decision variables increases quadratically with the problem horizon because of the feedback terms. 
Furthermore, integer variables that are used to encode obstacle avoidance constraints add to computational complexity.
In \cite{p:yin2021improvingMPPIwithCS}, the authors use unconstrained covariance steering to take sample trajectories from the low cost regions of the state space to enhance the performance and to avoid local minima. 
However, this method requires the terminal mean and covariance as design parameters, which can be hard to tune, and the safety constraints are not explicitly enforced.

\noindent\textbf{\text{Main Contributions:}} This paper presents a novel trajectory optimization algorithm for stochastic linear systems with a non-convex safe state space (CCSMPPI).
The proposed algorithm achieves high performance with safety guarantees by combining the standard MPPI with the CCS. 
In particular, MPPI is used to generate a reference trajectory, which is then used to generate a convex safe region, by solving an unconstrained stochastic trajectory optimization problem whereas CCS generates a control policy which minimizes the divergence from the reference trajectory while satisfying the safety constraints.
In this way, CCSMPPI endows the MPPI algorithm with robustness to stochastic disturbances by leveraging the framework of CCS. 

The other improvement of the CCSMPPI over the standard MPPI is the robustness and the safety guarantees against poorly designed cost functions and incorrect tuning of algorithm parameters.
The CCS procedure of the proposed algorithm filters the unsafe inputs that are computed by MPPI and corrects them by means of a feedback control law. 
This technique makes the cost function design and parameter tuning tasks less time-consuming. 

Finally, the practicality of CCSMPPI is demonstrated in two different trajectory optimization problems in which we compare our results with those obtained by using the standard MPPI \cite{p:williams2017informationMPPI} and tube-MPPI \cite{p:williams2018tubeMPPI}. 
It is shown in numerical experiments that our approach is superior to standard MPPI and tube-MPPI in terms of providing safety against both stochastic disturbances and poorly designed cost functions.

\section{Problem Statement and Preliminaries}\label{sec:problem-statement}
\subsection{Notation}
We denote by $\mathbb{R}^n$ the set of $n$-dimensional real vectors.
We use $\mathbb{E}\left[\cdot\right]$ and $\mathbb{P}(E)$ to denote the expectation functional and the probability of the random event $E$, respectively.
Given a vector $x$, its 2-norm is denoted by $\lVert x \rVert_2$ and
given a matrix $\mathbf{A} \in \mathbb{R}^{n \times m}$, we denote its Frobenius norm by $\lVert\mathbf{A} \rVert_{F}$ 
and its trace by $\operatorname{tr}(\mathbf{A})$. 
We use $\bm{0}$ and $I_n$ to denote the zero matrix 
and the $n \times n$ identity matrix, respectively.
We will denote the convex cone of $n\times n$ symmetric positive semi-definite (symmetric positive definite) matrices by $\mathbb{S}^{+}_n$ ($\mathbb{S}^{++}_n$). 
We write $\mathrm{bdiag}(A_1,$ $\dots, A_\ell)$ to denote the block diagonal matrix formed by the matrices $A_i$, $i\in \{1,\dots, \ell\}$. $\cup_{i=1}^{N} \mathcal{O}_i$ to denote the union of sets $\mathcal{O}_i$ indexed by $i \in \{1, \dots, N\}$.
We denote by $\mu_z$ and $\mathrm{var}_z$ the mean and the variance of a random vector $z$, respectively.
We write $ x \sim \mathcal{N}(\mu, \Sigma)$ to represent the random variable $x$ is normally distributed with mean $\mu$ covariance $\Sigma$.

\subsection{Problem Statement}
The motion of the agent is described by a discrete-time stochastic linear system:
\begin{align}
    x_{k+1} = A_k x_k + B_k u_k + w_k, \label{eq:linear-dynamics}
\end{align}
where $A_k \in \mathbb{R}^{4 \times 4}$, $B_k \in \mathbb{R}^{4 \times 2},$ $x_k \in \mathcal{X} \subseteq \mathbb{R}^{4}$ is the state of the agent which is decomposed as $x_k = [p_k^\mathrm{T}, v_k^\mathrm{T}]^{\mathrm{T}}$, 
where $p_k = [p^x_k, p^y_k]^{\mathrm{T}} \in \mathbb{R}^{2}$ is the position and $v_k = [v^x_k, v^y_k]^{\mathrm{T}} \in \mathbb{R}^{2}$ is the velocity,
$u_k = [u_k^x, u_k^y]^{\mathrm{T}} \in \mathbb{R}^2$ is the control input and 
$w_k \in \mathbb{R}^{4}$ is the disturbance. 
We assume that $w_k \sim \mathcal{N}(\bf{0}, \mathbf{W}_k)$ and $\mathbb{E}[w_k w_l^{\mathrm{T}}] = \bm{0}$ for all $k \neq l$.

The choice of the objective function will be determined by the particular application. However, in this paper, we use the objective function utilized in the formulation of the information-theoretic MPPI \cite{p:williams2017informationMPPI}, which is defined as follows:
{\small{
\begin{align}\label{eq:objective-function}
    \mathcal{L}(X^{N}, U^{N-1}) = \Phi(x_{N}) + \sum_{k=0}^{N-1} \left( q(x_k) + \lambda u_k^{\mathrm{T}} R_k u_k \right),
\end{align}
}}
\noindent where $X^N = \{ x_0, x_1, \dots, x_N \}$ is the state sequence, $U^{N-1} = \{ u_0, u_1, \dots, u_{N-1} \}$ is the control input sequence, $q: \mathbb{R}^{n} \rightarrow \mathbb{R}$ is the state-dependent term of the running cost function, $\Phi : \mathbb{R}^{2n} \rightarrow \mathbb{R}$ is the terminal cost function, $R_k \in \mathbb{S}_n^{+}$ and $\lambda \geq 0$. While the term $\lambda u_k^{\mathrm{T}} R_k u_k$ penalizes the control input, $q(x_k)$ and $\Phi(x_N)$ are task dependent and each of them can either be a smooth function as in \cite{p:williams2017informationMPPI} or a sum of indicator functions for obstacle avoidance as in \cite{p:williams2018tubeMPPI}.
 
We assume that the position space is populated by  $N_{\mathrm{obs}}$ obstacles.
The $i$-th obstacle is parametrized by its position, $s_i \in \mathbb{R}^{n}$, and its radius, $r_i \geq 0 $, 
and the region it occupies is denoted as $\mathcal{O}_i := \{ p \in \mathbb{R}^{n} |  \lVert x-s_i \rVert_{2} \leq r_i\}$.
The position space is defined as $\mathcal{X} \subseteq \mathbb{R}^{n}$ and the safe region is defined as $\mathcal{X}_{\mathrm{safe}} := \mathcal{X}\backslash\mathcal{O}$ where $\mathcal{O} = \cup_{j=1}^{N_{\mathrm{obs}}} \mathcal{O}_j $. Then, the safe trajectory optimization problem can be formally stated as follows:
\begin{subequations}\label{prob:prob1}
\begin{align}
    \underset{X^{N}, \; U^{N-1}}{\textrm{minimize}} & \quad \mathbb{E} [\mathcal{L}(X^{N}, U^{N-1})] \\
    \textrm{subject to}  & \quad x_0= \Bar{x}_0, ~~ \eqref{eq:linear-dynamics}  \nonumber  \\
    & \quad \mathbb{P}(x_k \in \mathcal{X}_{\mathrm{safe}}) \geq 1 - \mathrm{P}_{\mathrm{fail}}, ~~ \forall k \in \mathcal{I}^{N}
\end{align}
\end{subequations}
where $\mathcal{I}^{N} =  \{ 0, 1, \dots, N \}$, $\Bar{x}_{0} \in \mathbb{R}^{2n}$ and $\mathrm{P}_{\mathrm{fail}} \in (0,0.5]$ is the acceptable level of the probability of the failure. 

\noindent Since the system dynamics in \eqref{eq:linear-dynamics} includes the random noise term $w_k$, the state sequence $X^{k}$ is a random process thus the problem in \ref{prob:prob1} is a stochastic trajectory optimization problem.


\section{Brief Review of MPPI and TUBE-MPPI}\label{sec:MPPI}
MPPI is a sampling-based stochastic Model Predictive Control (MPC) algorithm \cite{p:williams2017informationMPPI}.
It works by taking  $K$ samples of control sequences from a Gaussian distribution, and finding the corresponding state trajectories and costs. 
Each sequence is then weighted by an exponential transform of its cost and the optimal control sequence is found as as weighted sum:
$\bm{u}^{\mathrm{MPPI}} = \frac{1}{\eta} \sum_{i = 1}^{K} \omega_{i} \bm{u}^{(i)}$, where $\eta = \sum_{i=1}^{K} \omega_i$
, $\omega_i = \exp{ -\frac{1}{\lambda} (C_i - \beta)},$
$\beta = \min_{i =1, \dots, K} C_i$, $\bm{u}^{(i)} =  \bm{v} + \bm{\epsilon}^{(i)}$ and $\bm{\epsilon}^{(i)} = [\epsilon^{(i)}_{0}, \dots, \epsilon^{(i)}_{T_{\mathrm{MPPI}}}]$, $\epsilon^{(i)}_{k} \sim \mathcal{N}(\bm{0}, \nu I)$ and $C_i$ is defined as follows:
\begin{align}
    & C_i = \Phi (x^{(i)}_{T_{\mathrm{MPPI}}}) + \sum_{k=0}^{T_{\mathrm{MPPI}}} \Big(q(x_k^{(i)}) + \frac{1}{2} v^{\mathrm{T}}_k R_k v_k \nonumber \\ 
    & \qquad + \frac{1}{2} \big(  v_k^{\mathrm{T}} R_k \epsilon^{(i)}_{k} +(1-\nu^{-1}) \epsilon_{k}^{(i)\mathrm{T}} R_k \epsilon_{k}^{(i)}\big) \Big),
\end{align}
$C_i$ is the total path cost induced by $i$th sample trajectory, $\bm{v} =  [v_0, \dots, v_{T_{\mathrm{MPPI}} - 1}]^{\mathrm{T}}$ is the optimal control sequence obtained by the previous iteration of the MPPI, $\bm{\epsilon}^{(i)}$ is the control sampling noise, $u^{(i)} $ is the $i$th control sequence sample and $\nu \geq 0$ is the control sampling covariance parameter.
At the next iteration, the previous optimal control sequence is shifted in time and used as the mean of the Gaussian distribution to sample controls from.

MPPI may not perform well when the system has a disturbance that causes the control sampling distribution to only sample high cost trajectories. 
Tube-MPPI~\cite{p:williams2018tubeMPPI} addresses this by running MPPI from  two starting states, a nominal and a real state. 
The real state is taken from the state of the real system as before whereas the nominal state is found by propagating the previous nominal state and control through the noise-free dynamics model without the disturbance term, $w_k$. 
Once the MPPI optimization step is done, the optimal trajectory of the real system is pushed to follow the optimal trajectory of the nominal system through the use of a feedback controller.
This procedure makes the trajectory of the real system converge to the trajectory of the nominal system which has a lower cost.

\section{Linear Covariance Steering Theory}\label{sec:CovSteering}
The main objective of the constrained covariance steering (CCS) problem is to steer the mean and the covariance of a stochastic linear system to desired values while minimizing the expected value of an objective function subject to state and/or input constraints  \cite{p:chen2015optimal, p:covarianceSkelton}. 
The general form of the discrete-time CCS problem can be formally stated as follows:
\begin{subequations}
\begin{align}\label{prob:genericCovSteerProb}
    \underset{\pi \in \Pi}{\textrm{minimize}} & ~~ J(\bm{x}, \bm{u}) = \mathbb{E}\left[\sum_{k=0}^{N-1} x_k^{\mathrm{T}} Q_k x_k + u_k^{\mathrm{T}} R_k u_k \right] \\ 
    \textrm{subject to} & ~~ x_{k+1} = A_k x_k + B_k u_k + w_k, \\
    & ~~ \mathbb{P}[x_k \in \mathcal{X}] \geq 1 - \epsilon_x, ~ \mathbb{P}[u_k \in \mathcal{U}] \geq 1 - \epsilon_u  \label{eq:generic-constr-chance}\\
    & ~~ x_0 \sim \mathcal{N}(\mu_0, \Sigma_0), \quad x_N \sim \mathcal{N}(\mu_\mathrm{d}, \Sigma_\mathrm{d}) \label{eq:generic-constr-covariance}
\end{align}
\end{subequations}
where $\Pi$ denotes the set of causal policies, $\mathcal{X} \subseteq \mathbb{R}^{n_x}$ and $\mathcal{U} \subseteq \mathbb{R}^{n_u}$ are arbitrary sets corresponding to state and input constraints.
The constraints in \eqref{eq:generic-constr-covariance} represents the initial and the desired state distributions.
The terminal covariance constraint in \eqref{eq:generic-constr-covariance} which can be dropped in our CCS formulation since it is not useful to specify a desired covariance for safety as long as constraints in \eqref{eq:generic-constr-chance} are satisfied.

Discrete-time formulations of the CCS problems can be cast as finite-dimensional deterministic optimization problems by restricting the class of admissible policies to those which admit the affine state history feedback \cite{p:BAKOLAS2018} or 
the disturbance feedback parametrization \cite{p:balci2021covariancedisturbance}.
In this work, we will be utilizing the latter according to which the control input at each discrete stage can be expressed as follows:
\begin{align}\label{eq:inputCS}
    v_k =
    \left \{
    \begin{matrix}
    \Bar{v}_{0} + H_0 (x_0 - \mu_0) &   k = 0, \\
    \Bar{v}_{k} + H_k (x_0 - \mu_0) + K_{k-1} w_{k-1} & k > 0,
    \end{matrix} \right. 
\end{align}
where $H_k, K_k \in \mathbb{R}^{n \times n}$ and $\Bar{v}_k \in \mathbb{R}^{n}$ for all $k  \in \{0, 1, \dots, N-1 \}$. 
This parametrization allows us to cast the covariance steering problem as a finite-dimensional (deterministic) optimization problem in terms of the following decision variables: $\{ \Bar{v}_k, H_k, K_k\}_{k = 0}^{N-1}$.
In order to do that, the decision variables are represented in a more compact form as follows:
$\bm{\mathcal{H}} := [H_0^{\mathrm{T}}, H_1^{\mathrm{T}}, \dots, H_{N-1}^{\mathrm{T}}]^{\mathrm{T}}$, 
$\mathbf{K} = \mathrm{bdiag}(K_0, K_1, \dots, K_{N-2}, \bm{0}^{n \times n})$, 
$\bm{\mathcal{K}}:= \left[ \begin{smallmatrix} \bm{0} & \bm{0} \\ \mathbf{K} & \bm{0} \end{smallmatrix} \right]$ and $\bm{\Bar{v}} = [\Bar{v}_0^{\mathrm{T}}, \Bar{v}^{\mathrm{T}}_1, \dots, \Bar{v}^{\mathrm{T}}_{N-1}]^{\mathrm{T}}$.

Now, let $\bm{x}$, $\bm{v}$ and $\bm{w}$ be concatenated state, input and noise vectors. Then, it follows from \eqref{eq:linear-dynamics} that
\begin{subequations}
\begin{align}
    \bm{x} &:= \bm{\Gamma} x_0 + \Gu \bm{v} + \Gw \bm{w}, \label{eq:concatenated-dynamics} \\
    \bm{v} &:= \bm{\Bar{v}} + \bm{\mathcal{H}} \Tilde{x}_{0} + \bm{\cK} \bm{w}, \label{eq:concatenated-input}
\end{align}
where $\Tilde{x}_0 := x_0 - \mu_0$. Equations \eqref{eq:concatenated-dynamics}-\eqref{eq:concatenated-input} are derived from \eqref{eq:linear-dynamics} and \eqref{eq:inputCS} respectively. The reader can refer to \cite{p:BAKOLAS2018,p:balci2021covariancedisturbance} for the details of the previous derivations.
\end{subequations}

We can compute the mean of the vectors $\bm{x}$ and $\bm{u}$ by taking the expectation of both sides of \eqref{eq:concatenated-dynamics} and \eqref{eq:concatenated-input}. Then, we compute the deviation of concatenated vectors $\bm{\Tilde{x}} = \bm{x} - \mu_{\bm{x}}$ and $\bm{\Tilde{v}} := \bm{v} - \mu_{\bm{v}}$. Finally, we compute the variances  $\mathrm{var}_{\bm{x}} := \mathbb{E}[\bm{\Tilde{x}} \bm{\Tilde{x}}^{\mathrm{T}}]$ and $\mathrm{var}_{\bm{v}} := \mathbb{E}[\bm{\Tilde{v}} \bm{\Tilde{v}}^{\mathrm{T}}]$ as follows: 
\begin{subequations}\label{eq:mean-cov-all}
\begin{align}
    \mu_{\bm{x}}&=\bm{\Gamma} \mu_{0} + \Gu \bm{\Bar{u}}, \label{eq:meanX}\\
    \mathrm{var}_{\bm{x}} &= (\bm{\Gamma} + \Gu \bm{\mathcal{H}})  \Sigma_0 (\bm{\Gamma} + \Gu \bm{\mathcal{H}})\t \nonumber \\ 
    & ~ ~ + (\Gw + \Gu \bm{\mathcal{K}}) \mathbf{W} (\Gw + \Gu \bm{\mathcal{K}})\t, \label{eq:covX}\\
    \mu_{\bm{v}} & = \bm{\Bar{v}}, \label{eq:meanV}\\ 
    \mathrm{var}_{\bm{v}} & = \bm{\mathcal{H}} \Sigma_0 \bm{\mathcal{H}}\t + \bm{\mathcal{K}} \mathbf{W} \bm{\mathcal{K}}\t \label{eq:covV}, 
\end{align}
\end{subequations}
where $\mu_{0} = \mathbb{E}[x_0]$, $\Sigma_0 = \mathbb{E}[\Tilde{x}_{0} \Tilde{x}_{0}^{\mathrm{T}}]$ and $\mathbf{W} = \operatorname{bdiag}(\mathbf{W}_{0}, \dots, \mathbf{W}_{N-1} )$.
Furthermore, the state and the control at the discrete stage $k$ can be recovered from the concatenated state and input vectors as $x_k = F^x_k \bm{x}$ and $v_k = F^{v}_{k}\bm{v}$, where $F^{x}_k$ and $F^{v}_{k}$ denote the block matrices whose $k$th block is equal to the identity matrix and the other blocks are equal to zero. Thus, the mean and the covariance of $x_k$ and $v_k$ are given by $\mu_{x_k} = F^{x}_{k} \mu_{\bm{x}}$, $\mathrm{var}_{x_k}  = F^{x}_{k} \mathrm{var}_{\bm{x}} F_k^{x\mathrm{T}}$, $\mu_{v_k} = F^{v}_{k} \mu_{\bm{v}}$, $ \mathrm{var}_{v_k} = F^{v}_{k} \mathrm{var}_{\bm{v}} F_k^{u\mathrm{T}}$.

So, it follows from \eqref{eq:mean-cov-all} that the mean of $x_k$ and $v_k$ can be expressed as affine functions of the decision variable $\bm{\Bar{v}}$ whereas the covariance matrices of $x_k$ and $v_k$ can be expressed as convex quadratic functions of the decision variables $\bm{\mathcal{H}}$ and $\bm{\mathcal{K}}$. 
This allows us to cast various forms of the (constrained) covariance steering problem as convex optimization problems which can be solved with highly efficient solvers. 
The reader can refer to \cite{p:BAKOLAS2018, p:balci2021covariancedisturbance, p:bakolas2017covIncompleteState, p:pilipovsky2021CovIteartiveRisk} for more details.

\section{Main Algorithm}\label{sec:MainAlgorithm}
The main components of the algorithm are the MPPI controller, the half-space generator, and the Constrained Covariance Steering module. 
The MPPI controller solves the unconstrained stochastic trajectory optimization problem and returns a state and an input sequence of length $T_{\mathrm{MPPI}}$. 
The state sequence generated by the MPPI module is used to generate the half-space constraints. 
The state and the input sequences and half-space constraints are used in the CCS module to solve for a policy that is guaranteed to be safe with high probability. 


\begin{algorithm}
    \SetKwInOut{Input}{Require} \SetKwInOut{Output}{output}
    \SetKwFunction{halfspace}{HalfSpaceGen}
    \SetKwFunction{mppi}{MPPI}
    \SetKwFunction{covsteer}{LinCovSteer}
    \SetKwFunction{setnom}{SetNomState}
    \SetKwFunction{pub}{Publish}
    \Input{$T_{\mathrm{max}}$, $T_{\mathrm{CS}}$, $T_{\mathrm{MPPI}}$, $\mathcal{M}_{\mathrm{obs}}, \{A_k, B_k, \mathbf{W}_k \}_{k=0, \dots, T_{\mathrm{max}}}, \sigma_{\mathrm{max}}$}
    
    $\Bar{x}_{0} \leftarrow x_{0}$\;
    $\Sigma_{k} \leftarrow \bm{0}$\;
    \For{$k \in \{0, 1, \dots, T_{\mathrm{max}}\}$}{
    $X^{T_{\mathrm{MPPI}}}, U^{T_{\mathrm{MPPI}}} \leftarrow $ MPPI({$\Bar{x}_k$})\;
    $ \mathcal{S}_{k}^{\mathrm{obs}}\leftarrow $ HSGen($X^{T_{MPPI}},  \mathcal{M}_{\mathrm{obs}}$) \;
    
    $ \mu_0 \leftarrow \Bar{x}_{k} $; $\Sigma_0 \leftarrow \Sigma_k$ \;
    
    $\bm{\Bar{v}}, \bm{\mathcal{H}}, \bm{\mathcal{K}} \leftarrow$ CCS({$\mu_0, \Sigma_0, X^{T_{\mathrm{MPPI}}}, U^{T_{\mathrm{MPPI}}}, \mathcal{S}_{k}^{\mathrm{obs}}$})\; 
    $\Bar{u}_{k} \leftarrow \Bar{v}_0$, $L_k \leftarrow H_0$ \;
    
    $u_k \leftarrow \Bar{u}_{k} + L_k (x_k - \Bar{x}_{k})$ \;
    SendToActuators($u_k$)\;
    $\Bar{x}_{k+1} \leftarrow A_k \Bar{x}_k + B_k \Bar{u}_{k}$ \;
    $\Sigma_{k+1} \leftarrow (A_k + B_k L_k) \Sigma_k (A_k + B_k L_k)^{\mathrm{T}} + \mathbf{W_k} $ \;
    
    \If{$\lambda_{\mathrm{max}}(\Sigma) > \sigma_{\mathrm{max}}$}
        {$\Bar{x}_{k+1} \leftarrow x_k$;
        $\Sigma_{k+1} \leftarrow \bm{0}$; }
    
    }
\caption{CCSMPPI}\label{al:CSMPPI}
\end{algorithm}
\subsection{MPPI}\label{subsec:MPPI}
In this paper, we follow the procedures described in \cite{p:williams2017informationMPPI} to use MPPI algorithm. 
The MPPI module requires system dynamics and initial state $x_0$ to sample trajectories. 
The MPPI horizon $T_{\mathrm{MPPI}}$, the input sampling covariance parameter $\nu$, and the number of trajectory samples $K$ are required as algorithm parameters. 

\subsection{Half-space Generation}\label{subsec:halfspacegen}
Safe half-spaces are generated by the ``HSGen'' procedure.
It takes the first $T_{\mathrm{CS}}$ states of the sequence $X^{T_{\mathrm{MPPI}}}$ and projects the position vectors $p_\ell$ onto each obstacle $\mathcal{O}_{j}$. 
Then, it computes the supporting hyperplane $\mathcal{H}_{\ell, j} := \{ p \in \mathbb{R}^{2} ~|~ a_{\ell,j}^{\mathrm{T}} p - b_{\ell,j} = 0 \},$
at the point of projection such that $a_{\ell, j}^{\mathrm{T}} p_\ell - b_{\ell, j} \geq 0 \Rightarrow p_\ell \notin \mathcal{O}_{j}.$
The projection of position $p_\ell$ at time $\ell$ onto the obstacle $\mathcal{O}_j$ is denoted as $z_{\ell, j}$ and is defined as $z_{\ell,j} := s_j + h_{\ell, j} r_j$,
where $h_{\ell,j}:= (p_\ell - s_j)/\lVert p_\ell - s_j \rVert_2$ is the unit normal vector to obstacle $\mathcal{O}_j$ at the point $z_{\ell, j}$
towards $p_\ell$. 
We set $ a_{\ell,j} = h_{\ell,j}$ and, $b_{\ell,j} = a_{\ell,j}^{\mathrm{T}} z_{\ell, j}$.
So, we can express $a_{\ell,j}$ and $b_{\ell,j}$ in terms of $p_l, s_j$ and $r_j$ as follows:
{\small{
\begin{align}\label{eq:halfspace}
    a_{\ell,j} &= \frac{(p_\ell - s_j)}{\lVert p_\ell - s_j \rVert_2}, & b_{\ell,j} & = \frac{(p_\ell - s_j)^{\mathrm{T}} s_{j} }{\lVert p_\ell - s_j \rVert_2} + r_j.
\end{align}
}}
The halfspace generation process is repeated for each obstacle $\mathcal{O}_j$ and every time step $\ell$. The halfspace parameters are gathered in the set of tuples $\mathcal{S}_k^{\mathrm{obs}} = \{ (a_{\ell,j}, b_{\ell,j}) \}_{\ell=0,\dots,T_{MC}}^{ j=1,\dots,N_{\mathrm{obs}}}  $ to be used in Constrained Covariance Steering.
The condition $p_{\ell} \notin \mathcal{O}_j$ is not necessary for the half-space generation procedure.
Even if  $p_\ell \in \mathcal{O}_j$, the procedure described by equations in \eqref{eq:halfspace} generates a half-space $\mathcal{H}_{\ell, j}$ such that $a_{\ell,j}^{\mathrm{T}} p - b_{\ell,j} \leq 0$ holds for all $p \in \mathcal{O}_{j}$.

\subsection{Constrained Covariance Steering}\label{subsec:covsteer}
The goal of the Constrained Covariance Steering module is to minimize the deviation of the actual state and control sequence from the reference state and control sequences which are computed by the MPPI algorithm while satisfying the safety constraints. This problem can be formally stated as the following stochastic optimal control problem: 

\begin{subequations}\label{prob:problem2}
\begin{align}
    \underset{\pi \in \Pi}{\textrm{minimize}} & \quad J(\pi):= \mathbb{E} \Big[\sum_{\ell=0}^{T_{CS}-1} \delta x_{\ell}^{\mathrm{T}} Q_{\ell} \delta x_{\ell}  + \delta u_{\ell}^{\mathrm{T}} R_{\ell} \delta u_{\ell} \nonumber   \\
    &  \qquad \qquad \quad  + \delta x_{T_{CS}}^{\mathrm{T}} Q_{T_{CS}} \delta x_{T_{CS}} \Big]  \label{prob2:obj}\\
    \textrm{subject to} & \quad x_{\ell+1} = A_\ell x_\ell + B_\ell u_\ell + w_\ell,  \quad \forall \ell \in \mathcal{I}_{t} \label{prob2:dynamics}\\
    & \quad u_\ell = \pi(x_{0}, \dots, x_{\ell}), \qquad\qquad\; \forall \ell \in \mathcal{I}_{\mathrm{t}}\\
    &  \quad \mathbb{P} \left[ a_{\ell, j}^{\mathrm{T}} p_\ell - b_{\ell,j} \geq 0 \right] \geq 1 - P_{\mathrm{fail}},  \quad \nonumber \\
    & \qquad \qquad\qquad\qquad\qquad\qquad \forall \{\ell, j\} \in \mathcal{I} \label{prob2:constrsafe}
\end{align}
\end{subequations}
\noindent where $\Pi$ denotes the set of all admissible control policies, $\delta x_{\ell} = x_{\ell} - x_{\ell}^{\mathrm{MPPI}}$, $\delta u_{\ell} = u_{\ell} - u_{\ell}^{\mathrm{MPPI}}$, $\mathcal{I}_{\mathrm{t}} := \{0, \dots, T_{CS}\}$, $\mathcal{I}_{\mathrm{o}} := \{1, \dots, N_{\mathrm{obs}}\}$, $\mathcal{I} = \mathcal{I}_{\mathrm{t}} \times \mathcal{I}_\mathrm{o}$.
The stochastic optimal control problem defined in 
\eqref{prob:problem2} 
can be cast as a deterministic optimization problem by fixing the policy as in \eqref{eq:inputCS} 
as explained in Section \ref{sec:CovSteering}. 
The resulting finite-dimensional deterministic optimization problem is given by:
\begin{subequations}\label{prob:problem3}
\begin{align}
    \underset{\bm{\Bar{u}}, \bm{\mathcal{H}}, \bm{\mathcal{K}} }{\textrm{minimize}} & \quad \mathcal{J}(\bm{\Bar{u}}, \bm{\mathcal{H}}, \bm{\mathcal{K}}):=\delta \bm{\Bar{x}}^{\mathrm{T}} \bm{Q} \delta \bm{\Bar{x}} + \delta \bm{\Bar{u}}^{\mathrm{T}} \bm{\mathcal{R}} \delta \bm{\Bar{u}} \nonumber\\
    & \qquad \qquad + \operatorname{tr}(\bm{\mathcal{Q}} \mathrm{var}_{\bm{x}} ) + \operatorname{tr}(\bm{\mathcal{R}} \mathrm{var}_{\bm{u}}) \label{prob3:obj} \\
    \textrm{subject to} & \quad a_{\ell,j}^{\mathrm{T}} \mathbf{P}_{\ell} \mu_{\bm{x}} - b_{\ell,j} \geq \nonumber \\ 
    & \qquad \qquad \alpha \lVert \bm{\zeta}^{\mathrm{T}} \mathbf{P}_{\ell}^{\mathrm{T}} a_{\ell,j} \rVert_2, ~~ \forall \{\ell, j\} \in \mathcal{I} \label{prob3:constr1}\\
    & \quad \bm{\zeta} = \left[ (\mathbf{G}_0 + \mathbf{G_u} \bm{\mathcal{H}}) ~ (\mathbf{G_w} + \mathbf{G_w} \bm{\mathcal{K}}) \right]\mathbf{R} \label{prob3:constr2}
\end{align}
\end{subequations}
\noindent where $\mathbf{R}\mathbf{R}^{\mathrm{T}} = \operatorname{bdiag}(\Sigma_0, \mathbf{W})$, $\bm{x}$ and $\bm{u}$ are defined as in Section \ref{sec:CovSteering}, $\delta \bm{\Bar{x}} = \mu_{\bm{x}} - \bm{x}^{\mathrm{MPPI}}$, $\delta \bm{\Bar{u}} = \mu_{\bm{u}} - \bm{u}^{\mathrm{MPPI}}$, $\bm{\mathcal{Q}}=\operatorname{bdiag}(Q_0, \dots, Q_{T_{CS}})$, $\bm{\mathcal{R}}=\operatorname{bdiag}(R_0, \dots, R_{T_{CS}-1})$. 
$\mathbf{P}_{\ell}$ is defined such that $p_\ell = \mathbf{P}_{\ell} \bm{x}$
and $\alpha = \varphi^{-1}(1 - P_{\mathrm{fail}})$ where $\varphi$ is the cumulative density function of normal random variable with zero mean and unit variance.
Finally, $\mu_{\bm{x}}$, $\mu_{\bm{u}}$, $\mathrm{var}_{\bm{x}}$ and $\mathrm{var}_{\bm{u}}$ are defined as in \eqref{eq:meanX}-\eqref{eq:covV}. In addition, we observe that $\mathrm{var}_{\bm{x}} = \bm{\zeta} \bm{\zeta}^{\mathrm{T}}$ where $\bm{\zeta}$ is given by \eqref{prob3:constr2}.

To see the equivalence of optimization problems in \eqref{prob:problem3} and \eqref{prob:problem2}, 
first, observe that the objective function in \eqref{prob:problem2} can be written as $\mathbb{E}[\delta \bm{x}^{\mathrm{T}} \bm{\mathcal{Q}} \delta \bm{x} + \delta \bm{u}^{\mathrm{T}} \bm{\mathcal{R}} \delta \bm{u} ]$ 
where $\delta \bm{x} = \bm{x} - \bm{x}^{\mathrm{MPPI}}$, $\delta \bm{u} = \bm{u} - \bm{u}^{\mathrm{MPPI}}$.
Using the cyclic permutation property of trace operator, the linearity of expectation and the equalities $\mathrm{var}_{\delta \bm{x}} = \mathrm{var}_{\bm{x}}$ and $\mathrm{var}_{\delta \bm{u}} = \mathrm{var}_{\bm{u}}$, it follows that the objective functions in \eqref{prob2:obj} and \eqref{prob3:obj} are equivalent.
Proposition \ref{prop:chanceconstraint} along with the expressions $\mu_{p_\ell} = \mathbf{P}_{\ell} \mu_x$ and $\mathrm{var}_{p_{\ell}} = \mathbf{P}_{\ell} \bm{\zeta} \bm{\zeta}^{\mathrm{T}} \mathbf{P}_{\ell}^{\mathrm{T}}$ show that \eqref{prob3:constr1}, \eqref{prob3:constr2} $\Leftrightarrow$ \eqref{prob2:constrsafe}.
\begin{proposition}\label{prop:chanceconstraint}
Let $ p \sim \mathcal{N}(\mu_{p}, \Sigma_{p})$, where $a, \mu_p \in \mathbb{R}^n$, $b \in \mathbb{R}$, $\Sigma_p \in \mathbb{S}_n^{+}$, and $P \in (0, 1/2]$. Then, $\mathbb{P}[ a^{\mathrm{T}} p - b \geq 0] \geq 1 - P$  if and only if $a^{\mathrm{T}} \mu_p - b \geq \varphi^{-1} (1-P) \lVert \mathbf{R}^{\mathrm{T}} a \rVert_2 $ where $\varphi : \mathbb{R} \rightarrow (0, 1) $ is the cumulative density function of normally distributed random variable with zero mean and unit variance and, finally $\mathbf{R}$ is such that $\mathbf{R} \mathbf{R}^{\mathrm{T}} = \Sigma_p$.
\end{proposition}

The problem in \eqref{prob:problem3} has a convex quadratic objective function and the constraint in \eqref{prob3:constr2} is affine. 
Also, the constraint in \eqref{prob3:constr1} corresponds to a second-order cone constraint since $\alpha = \varphi^{-1}(1-P_{\mathrm{fail}}) \geq 0$ for all $P_{\mathrm{fail}} \in (0, 0.5]$ \cite{b:boyd2004convex}. 
Thus, problem \ref{prob:problem3} can be solved for global optimal solution $( \bm{\Bar{v}}^{\star}, \bm{\mathcal{H}}^{\star}, \bm{\mathcal{K}}^{\star})$ using off-the-shelf solvers such as \cite{mosek2010mosek}. 




\subsection{Discussion}\label{subsec:discussionAlgorithm}
First, the algorithm is initialized by setting $\Bar{x}_{0} = x_0$ and $\Sigma_0 = \bm{0}$ where $\Bar{x}_0$ and $\Sigma_0$ represent the initial nominal state and initial covariance respectively. 
Then, using $\Bar{x}_k$ as the initial state, MPPI generates a pair of reference state and input sequences $(X^{T_{\mathrm{MPPI}}}, U^{T_{\mathrm{MPPI}}})$. 
The state sequence is used to generate a safe convex region over which the constraints \eqref{prob2:constrsafe} are satisfied based on the technique that will be described in Section \ref{subsec:halfspacegen}. 
Then, we formulate a corresponding CCS problem that seeks for a control policy in the form of \eqref{eq:inputCS}
that will guarantee the satisfaction of chance constraints while minimizing the deviation from reference state and control sequences. 
If the largest eigenvalue of the computed covariance $\Sigma_{k+1}$ exceeds a predetermined threshold $\sigma_{\mathrm{max}}$, then the nominal state $\Bar{x}_{k}$ is set equal to the real state $x_k$ and covariance $\Sigma_k$ is set to $\bm{0}$. 
Next, the nominal state $\Bar{x}_{k+1}$ and covariance matrix $\Sigma_{k+1}$ will be updated as described in lines 9-10 in Algorithm \ref{al:CSMPPI}.

It is worth mentioning that the CCS module uses the disturbance noise covariance in its formulation. However, this information is usually unknown in real-world scenarios. But, this can easily be handled by over-approximating the noise covariance, that is, by taking $\mathbf{W}_k \succeq \mathbf{W}_k^{\mathrm{real}}$, where $\mathbf{W}_{k}^{\mathrm{real}}$ is the actual noise covariance that is acting on the system.
This allows the CCS module to find a policy that satisfies the safety constraints. 
Furthermore, this condition implies that the state will stay in the safe region with probability greater than $1 - P_{\mathrm{fail}}$.
Although this approach may generate overly conservative policies, system identification techniques can be used to learn the actual noise covariance \cite{p:feng2014kalman} and hence reduce conservativeness.



The final component of our algorithm is the use of nominal dynamics, which is the noise-free version of the real dynamics. 
Since the CCS module assures that chance constraints are satisfied and it uses the nominal state as the initial mean state in its formulation, 
the safety margins of the nominal state will be greater than the real state.
Also, by computing the covariance in line 10 in Algorithm \ref{al:CSMPPI}, we compute the high probability region where the real state lies. Then, this covariance value is used as initial covariance in the CCS procedure in the next step to guarantee the satisfaction of the chance constraints. 


\section{Numerical Experiments}\label{sec:Numerical Experiments}
In our numerical experiments, we consider a double integrator dynamics described by \eqref{eq:linear-dynamics} with euler discretization scheme with 
$dt = 0.05$, and $\{w_k\}$ is taken to be a white noise process with $w_k \sim \mathcal{N}(\bm{0}, \mathbf{W}_k)$ where the noise covariance matrix $\mathbf{W}_k$ varies depending on different problem instances. 
We show the efficacy of our approach in two trajectory optimization problems: 
an obstacle avoidance problem and a path generation problem in a circular track.


\noindent \textbf{Obstacle Avoidance:} 
In the obstacle avoidance case, we compare the performances of CCSMPPI with tube-MPPI \cite{p:williams2018tubeMPPI} under high noise that is acting upon the system to show the robustness of our approach.
In the tube-MPPI formulation, an LQG tracking controller is used to track nominal state and input sequences.
In our experiments, the LQG cost function parameters $Q^{\mathrm{LQG}}_k, R^{\mathrm{LQG}}_k$ are chosen to be equal to the cost function parameters used in the CCS formulation given in \eqref{prob2:obj}. Also, the failure parameter $P_{\mathrm{fail}}$ that is shown in \eqref{prob2:constrsafe} is taken to be $0.01$.

We consider the running cost function $q_{\mathrm{hard}}(p_k) = \lVert p_k - p_{\mathrm{des}} \rVert_2^{2} +  5000  \sum_{j=1}^{N_{\mathrm{obs}}} \mathbb{I}_{\mathcal{O}_j}(p_k)$
and $\Phi_{\mathrm{hard}}(x_{T}) = 0$ where $\mathbb{I}_{\mathcal{O}_j} :  \mathbb{R}^{2} \rightarrow \{0, 1\}$ is the indicator function of set $\mathcal{O}_j$. 

The parameters of the MPPI algorithm used in these experiments 
are $T_{\mathrm{MPPI}}=40$, $K=100$, $\lambda=0.1$, $\nu = 0.1$ and $\epsilon_k \sim \mathcal{N} ({\bm{0}, 0.001 I})$. 
In addition, the problem horizon parameter $T_{\mathrm{max}} = 200$ and the noise covariance matrix $\mathbf{W}_k = \operatorname{bdiag}(0.0, 0.0, 5.0, 5.0)$. 
The state-dependent term of the running cost function was taken to be $q(x_k) = 10 \, q_\mathrm{hard}(x_k)$ and the desired final position $p_{\mathrm{des}} = [2.0, 10.0]^{\mathrm{T}}$.



\begin{figure}[ht]
    \centering
    \begin{subfigure}{0.48\linewidth}
    \includegraphics[width=\linewidth]{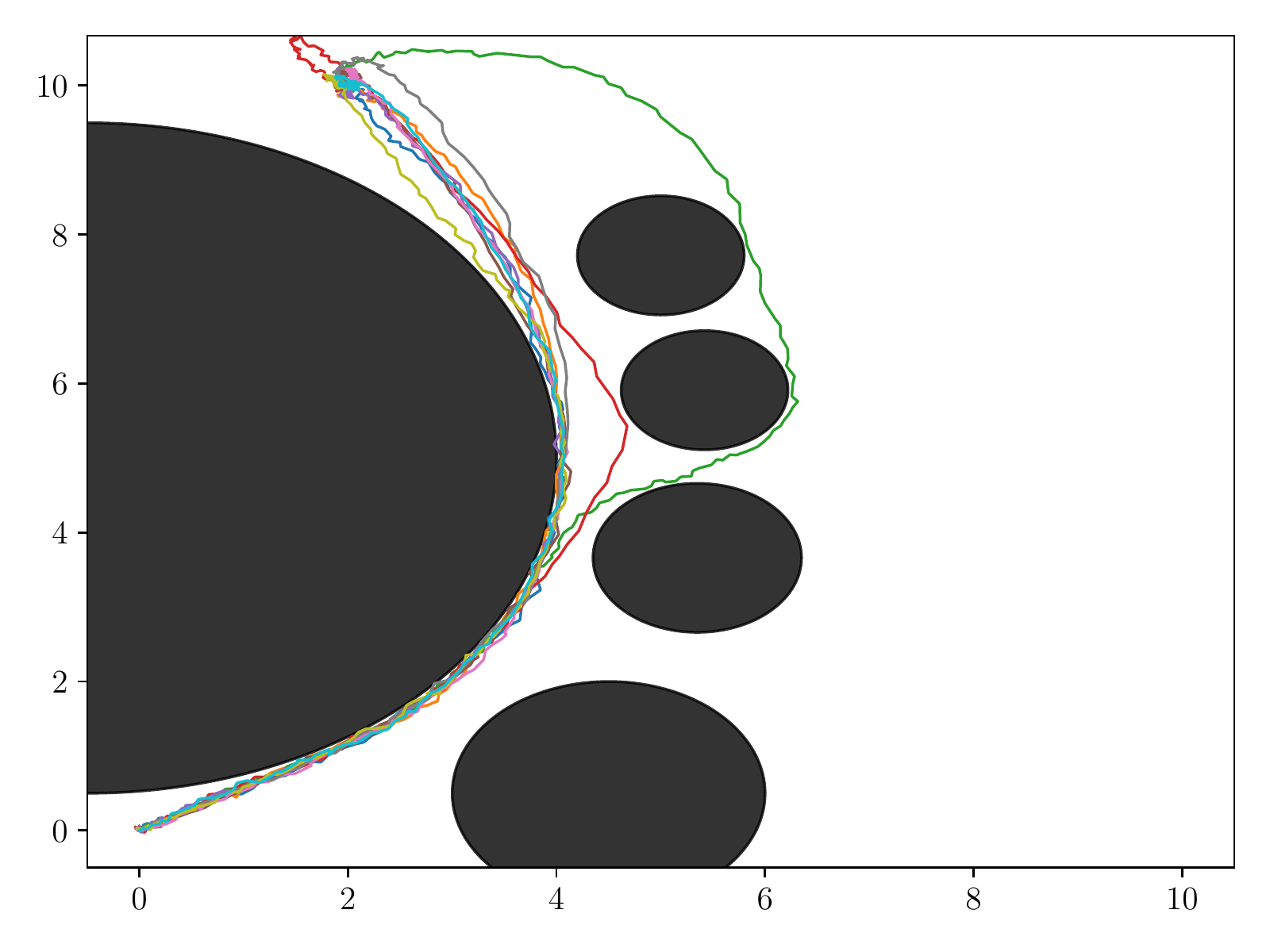}
    \caption{CCSMPPI}\label{fig:MChardCS}
    \end{subfigure}
    ~
    \begin{subfigure}{0.48\linewidth}
    \includegraphics[width=\linewidth]{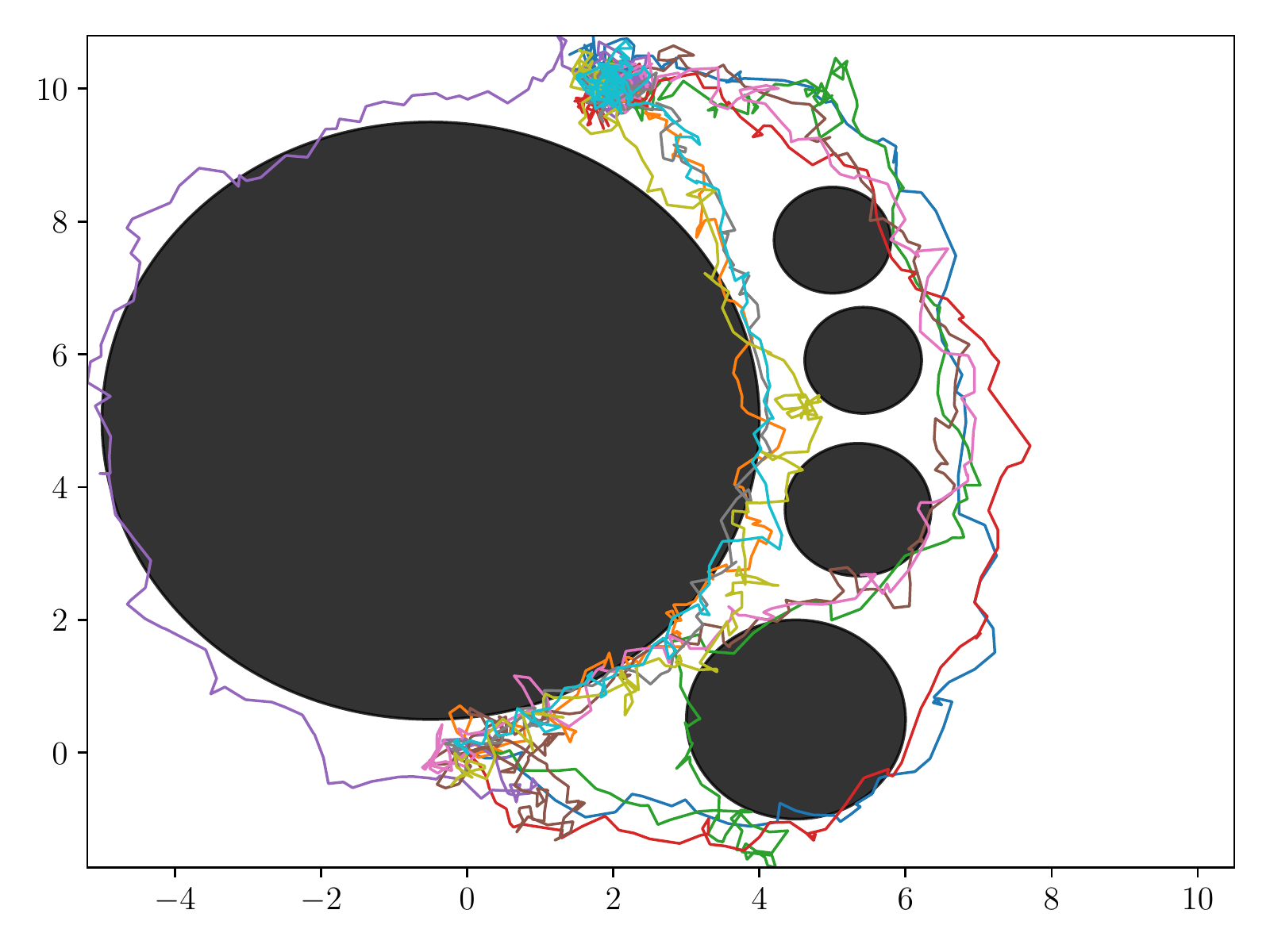}
    \caption{Tube-MPPI}\label{fig:MChardTUBE}
    \end{subfigure}
\end{figure}

\begin{figure}[ht]
    \centering
    \begin{subfigure}{0.35\linewidth}
    \includegraphics[width=\linewidth]{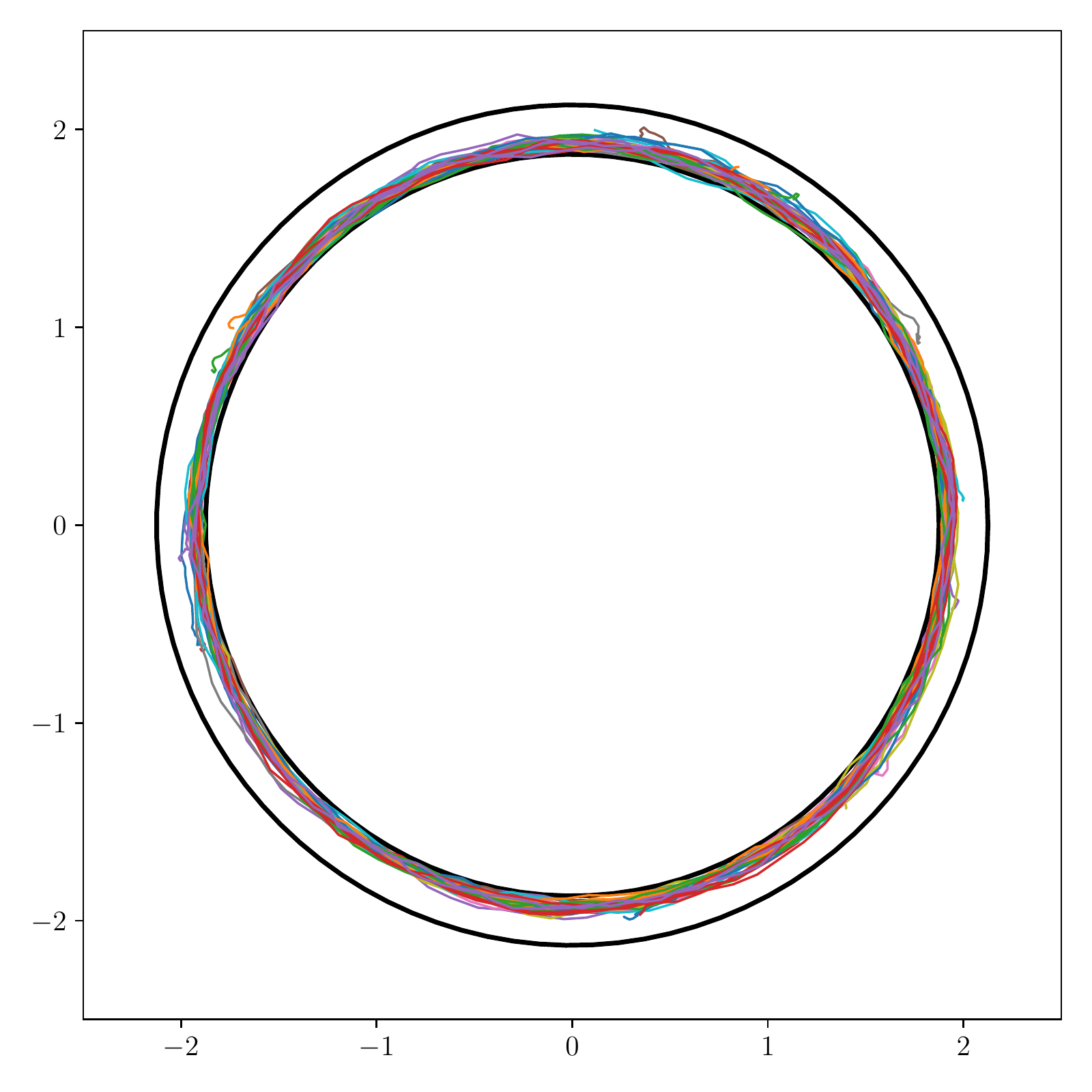}
    \subcaption{CCSMPPI - $q_{\mathrm{t,s}}$}\label{subfig:softCS}
    \end{subfigure}
    ~
    \begin{subfigure}{0.35\linewidth}
    \includegraphics[width=\linewidth]{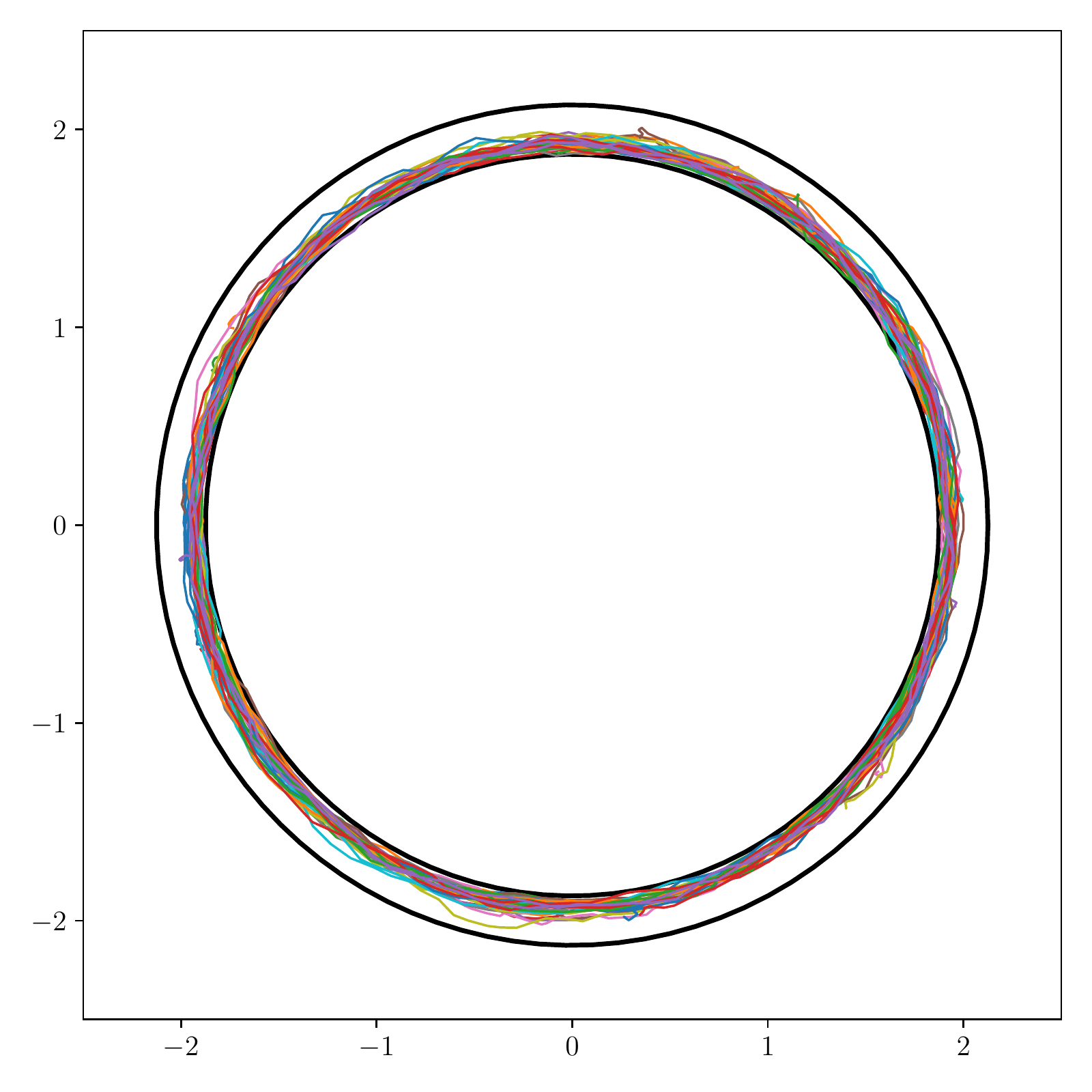}
    \subcaption{CCSMPPI - $q_{\mathrm{t,h}}$}\label{subfig:hardCS}
    \end{subfigure}
    ~
    \begin{subfigure}{0.35\linewidth}
    \includegraphics[width=\linewidth]{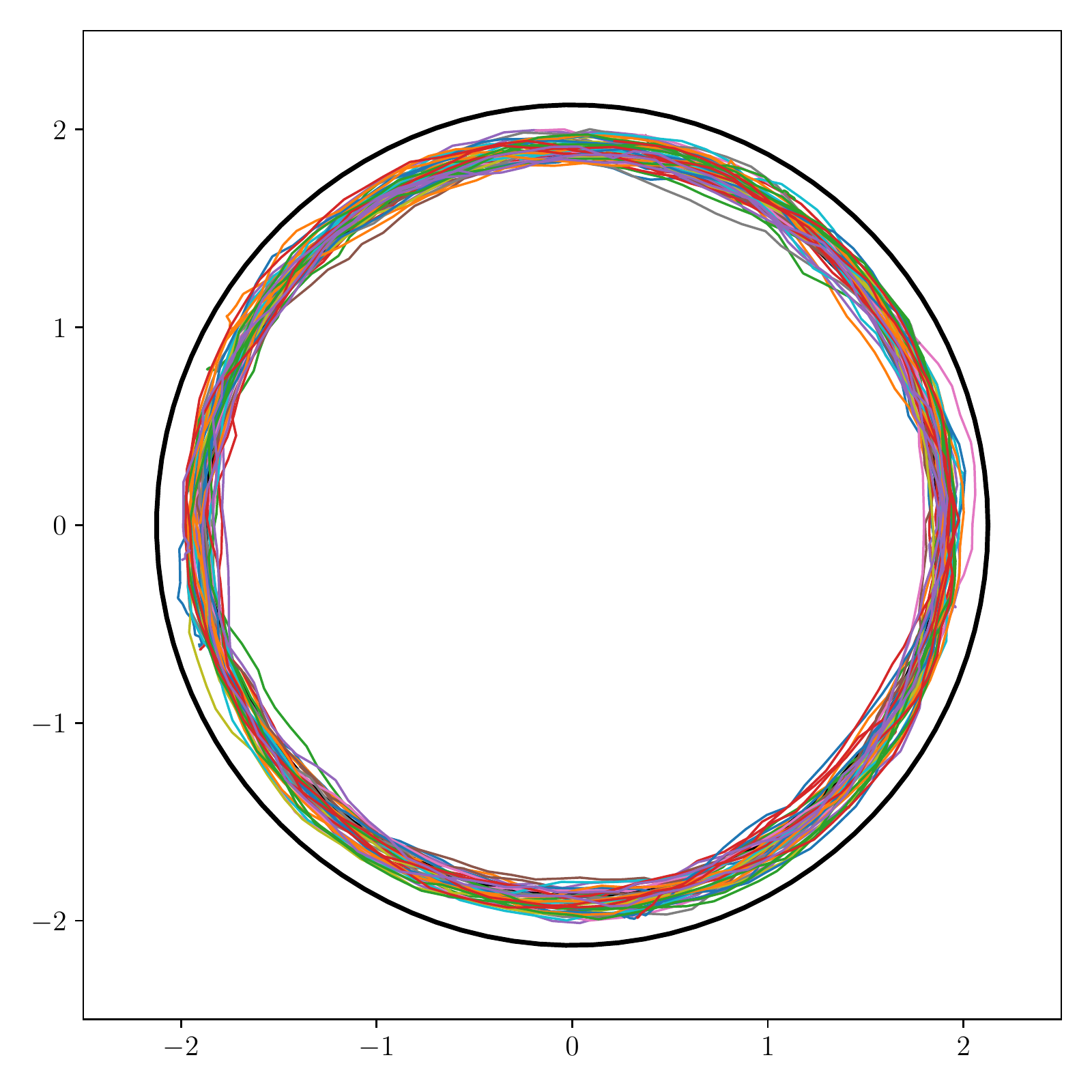}
    \subcaption{Tube-MPPI - $q_{\mathrm{t,s}}$}\label{subfig:softtube}
    \end{subfigure}
    ~
    \begin{subfigure}{0.35\linewidth}
    \includegraphics[width=\linewidth]{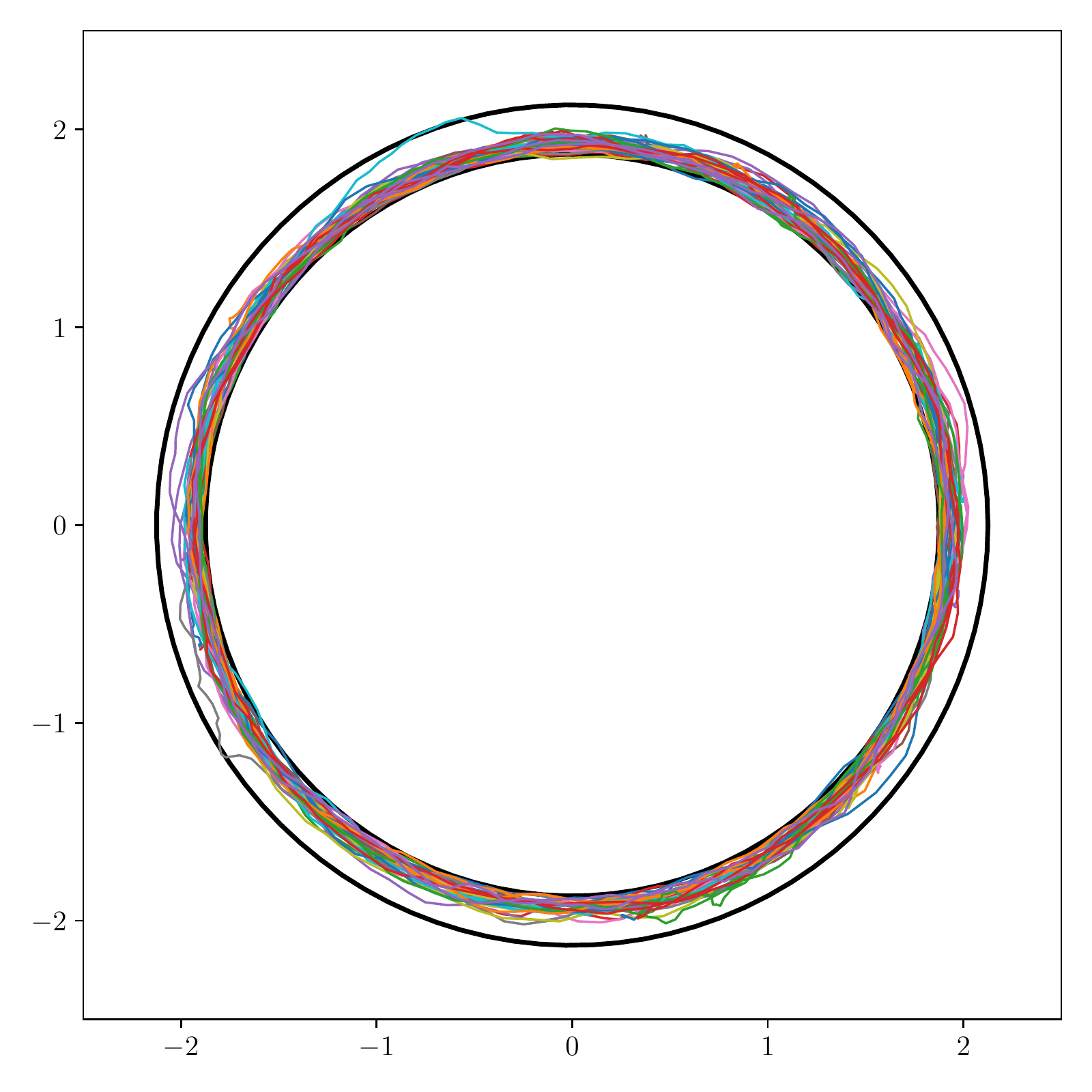}
    \subcaption{Tube-MPPI - $q_{\mathrm{t,h}}$}\label{subfig:hardtube}
    \end{subfigure}
    ~
    \begin{subfigure}{0.35\linewidth}
    \includegraphics[width=\linewidth]{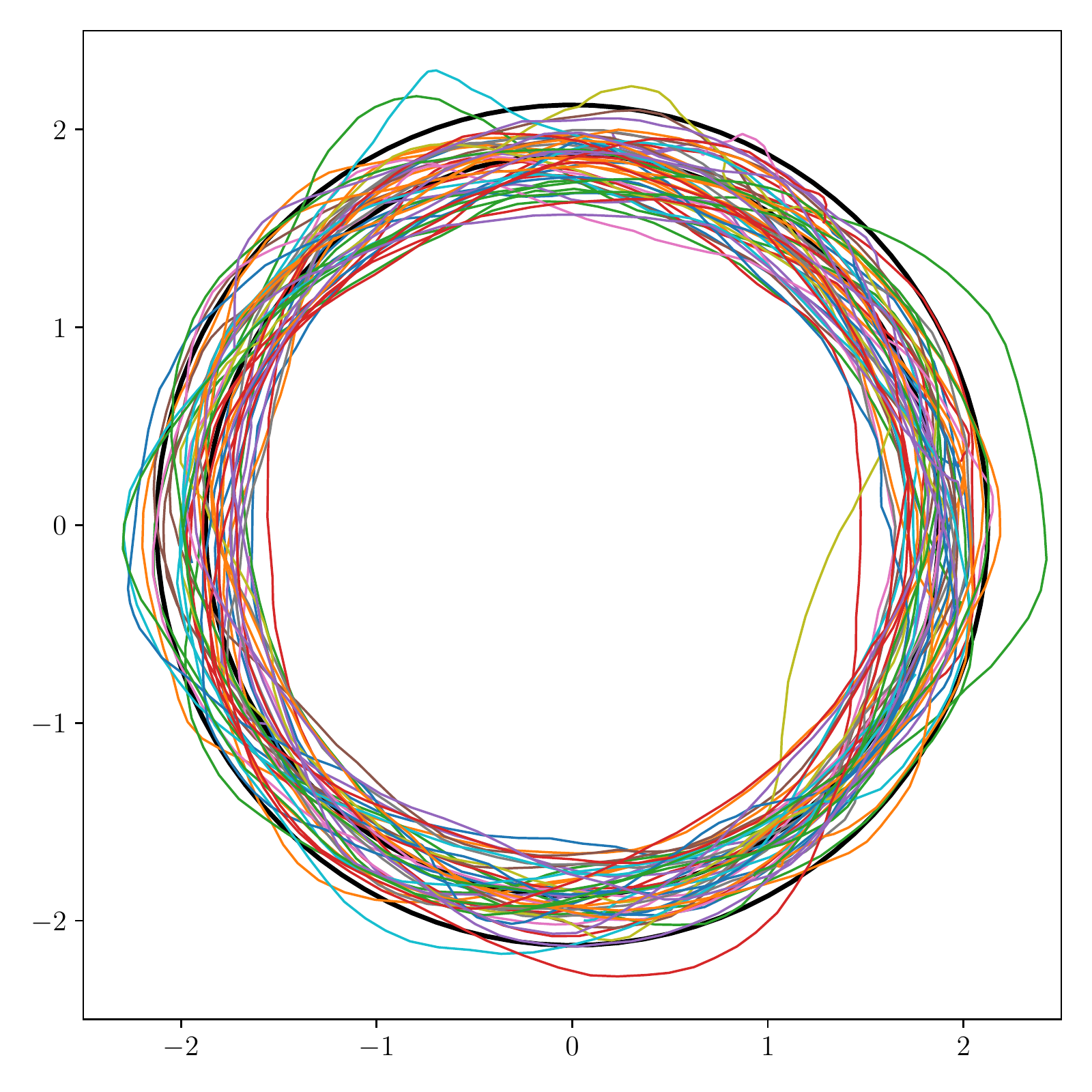}
    \subcaption{MPPI - $q_{\mathrm{t,s}}$}\label{subfig:softreg}
    \end{subfigure}
    ~
    \begin{subfigure}{0.35\linewidth}
    \includegraphics[width=\linewidth]{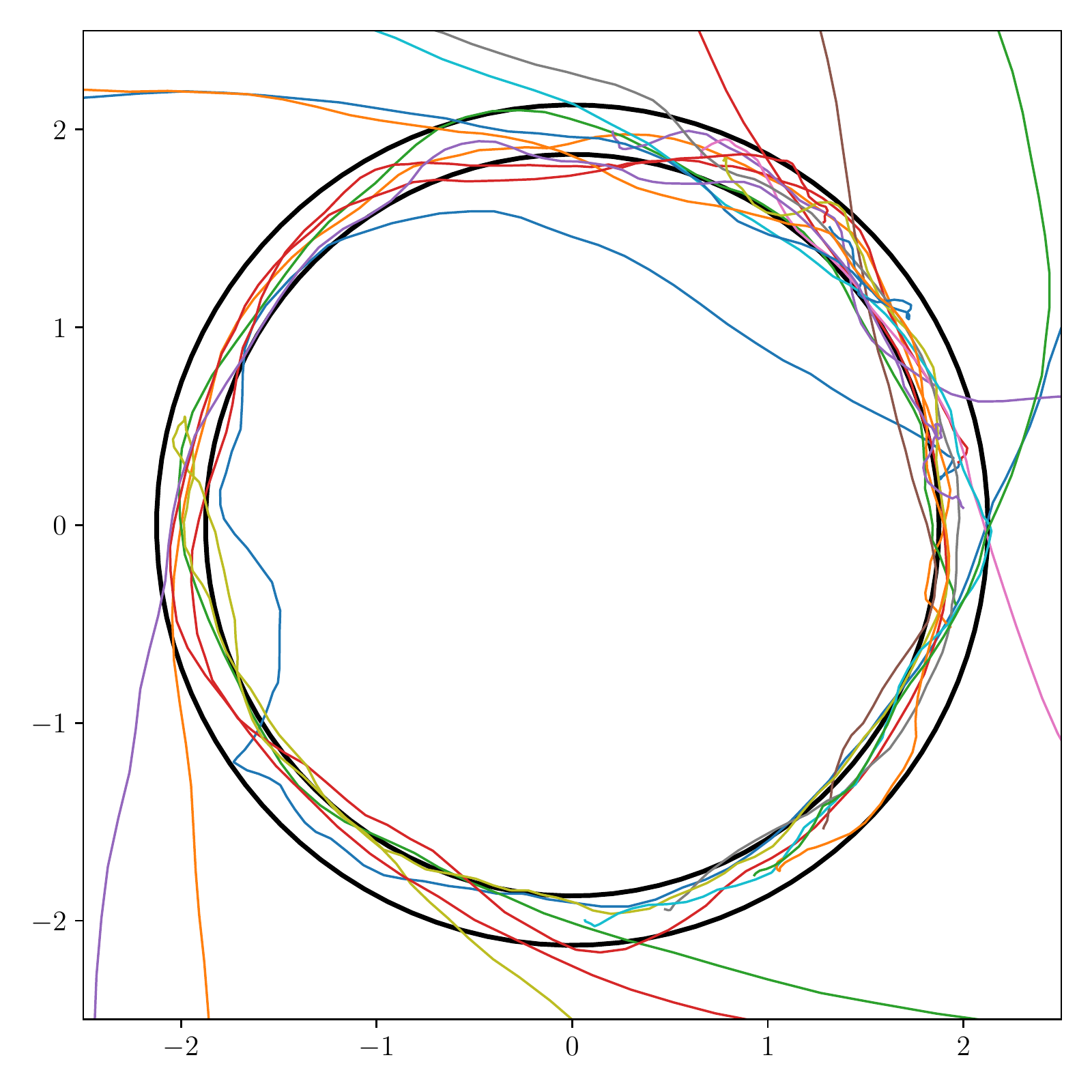}
    \subcaption{MPPI - $q_{\mathrm{t,h}}$}\label{subfig:hardreg}
    \end{subfigure}
    ~
    \caption{CCSMPPI, Tube-MPPI and standard MPPI in different scenarios 
    }\label{fig:comparison}
    \vspace{-0.7cm}
\end{figure}

\vspace{-0.3cm}
Figure \ref{fig:MChardCS} illustrates 10 randomly sampled trajectories induced by the CCSMPPI algorithm. 
Although the intensity of the noise that is acting upon the system is quite high compared to the sampling distribution parameter $\nu$, the CCSMPPI is successfully avoiding obstacles. 
Figure \ref{fig:MChardTUBE} shows 10 randomly sampled trajectories of the system running under the tube-MPPI algorithm. 
It can be seen that the agent reaches the goal position but fails to avoid obstacles even though 
the trajectories that lead to collisions are heavily penalized.
In this case, tube-MPPI fails to handle uncertain disturbances and generate safe control inputs.

\noindent \textbf{Circular Track:} In this scenario, the goal is to keep the position of the system in a circular track with inner and outer radius $R_{\mathrm{in}} , R_{\mathrm{out}} = R_{\mathrm{c}} \mp 0.125$ 
where $R_{\mathrm{c}}=2$ while maintaining a desired speed $v_{\mathrm{des}}$ in counter-clockwise direction.
In the numerical experiments, this goal is encoded as 2 different running cost functions $q_{\mathrm{t, s}}(x_k)$ and $q_{\mathrm{t,h}}(x_k)$:
\begin{subequations}
\begin{align}
    q_{\mathrm{t,s}} (x_k) &:= (\lVert v_k \rVert_2 - v_{\mathrm{des}})^2 + \lVert (p_k^x v^y_k - v^x_k p_k^y) - R_{\mathrm{c}} v_{\mathrm{des}}\rVert \nonumber \\
    & + 100 \left(\sqrt{p_k^{x2} + p_k^{y2}} - R_{\mathrm{c}} \right)^{2} \label{eq:qtracksoft}\\
    q_{\mathrm{t,h}} (x_k) &:= (\lVert v_k \rVert_2 - v_{\mathrm{des}})^2 + \lVert(p_k^x v^y_k - v^x_k p_k^y) - R_{\mathrm{c}} v_{\mathrm{des}}\rVert \nonumber \\ 
    & + 5000 \mathbb{I}_{\lnot C}(p_k) \label{eq:qtrackhard}
\end{align}
\end{subequations}
where $ C := \{ p \in \mathbb{R}^{2} \, | \, R_{\mathrm{c}} - 0.125 \leq \lVert p \rVert_2 \leq R_{\mathrm{c}} + 0.125  \}$,  $ \mathbb{I}_{\lnot C} : \mathbb{R}^2 \rightarrow \{0, 1\}$ is the indicator function such that $\mathbb{I}_{\lnot C} (p) = 1$ if $p \notin C$ and $\mathbb{I}_{\lnot C}(p) = 0 $ if $p \in C$. 

The safety criterion in this example is to stay within the circular track which is formally defined by set $C$. 
This condition is encoded in $q_{\mathrm{t,s}}(x_k)$ in \eqref{eq:qtracksoft} by penalizing the deviation of the position $p_k$ from the mid radius $R_{\mathrm{c}}$.
This choice makes $q_{t,s}(x_k)$ a smooth function.
On the other hand, the safety criterion is encoded by using indicator functions in \eqref{eq:qtrackhard} which is a more clear encoding of the safety constrained however non-smoothness of $q_{\mathrm{t, h}}(x_k)$ makes this problem harder to solve.

In Figure \ref{fig:comparison}, trajectories generated by CCSMPPI, Tube-MPPI and standard MPPI are shown from the top to the bottom. Figures \ref{subfig:softCS}, \ref{subfig:softtube}, \ref{subfig:softreg} show the results with running cost $q(x_k) = 100 \, q_{\mathrm{t,s}}(x_k)$ and Figures \ref{subfig:hardCS}, \ref{subfig:hardtube}, \ref{subfig:hardreg} show the results with $q(x_k) = 100 \, q_{\mathrm{t, h}}(x_k)$. 
The parameters of the MPPI are $T_{\mathrm{MPPI}}=20$, $K=200$, $\lambda=0.1$, $\nu = 1.0$, and $\epsilon \sim \mathcal{N}(\bm{0}, 0.001 I)$. 
In addition, the problem horizon $T_{\mathrm{max}} = 300$, $\mathbf{W_k} = \operatorname{bdiag}(0.001, 0.001, 1.0, 1.0)$ and $v_{\mathrm{des}} = 6.0$. 

When $q_{\mathrm{t,s}}(x_k)$ is used as a running cost function, the trajectories induced by tube-MPPI and standard MPPI fail to meet the safety criterion (shown in Figures \ref{subfig:softtube}, \ref{subfig:softreg}) due to the poor design of the cost function.
On the other hand, it is shown in Figure \ref{subfig:softCS} that the CCSMPPI manages to keep the position within the track. 
When the running cost is switched to $q_{\mathrm{t,h}}(x_k)$, the non-smoothness of the cost function causes standard MPPI to fail as shown in Figure \ref{subfig:hardreg}.
Even though tube-MPPI seems to keep the position within the circular track in Figure \ref{subfig:hardtube}, there are more violations of the safety constraints than CCSMPPI as is shown in Figure \ref{subfig:hardCS}.

In Table \ref{tab:Performance Comparison}, we compare the performance of CCSMPPI, tube-MPPI and standard MPPI by sampling $N_{\mathrm{sim}} =15 $ trajectories for both experiments \#1 and \#2. The running cost function is taken as $100 q_{t, s}(x_k)$ and $100 q_{t, h}(x_k)$ in experiments \#1 and \#2, respectively. 
Also, $T_{\mathrm{max}}$ is taken as $200$ and $300$ in experiments \#1 and \#2, respectively. 
In both experiments, $\mathbf{W_k}$ are chosen to be equal to $\operatorname{bdiag}(0.005, 0.005, 0.5, 0.5)$ and $v_{\mathrm{des}} = 6.0$. 
\textbf{Pr\textsubscript{fail}} represents the probability of failure and it is computed by dividing the number of trajectories that leave the circular track at least once ($N_{\mathrm{fail}}$) by the total number of trajectories $N_{\mathrm{sim}}$. 

It can be seen from the results of experiment \#1 in Table \ref{tab:Performance Comparison} that standard MPPI performs better in terms of minimizing the cost than both tube-MPPI and CCSMPPI and reaches higher speeds. 
However, this is due to the poor design of the cost function, and the fact that the control inputs that are corrected by CCS module to guarantee safety are not optimal with respect to the used cost function.
When $q_{\mathrm{t, h}}(x_k)$ is used as the running cost in experiment \#2, standard MPPI performs worse than both tube-MPPI and CCSMPPI due to the presence of random noise $w_k$.
In these experiments, the safety of the trajectory is the first priority, as encoded in the running cost $q_{\mathrm{t,h}}(x_k)$.
Although tube-MPPI reaches higher speeds, it fails to reach the safety levels of CCSMPPI. 

It should be highlighted that the probability of violating the constraint in \eqref{prob2:constrsafe} at every time step $k$ is less than $P_{\mathrm{fail}} = 0.01$ but still greater than 0.
This means that as $T_{\mathrm{max}} \rightarrow \infty$, the failure probability approaches 1. 
This is the reason why \textbf{Pr\textsubscript{fail}} is non-zero for CCSMPPI in both experiments. 
\textbf{Pr\textsubscript{fail}} can be reduced by lowering the safety threshold $P_{\mathrm{fail}}$, however it is not possible to make it 0 since $w_k$ is assumed to be normally distributed and thus unbounded. 

\begin{table}[ht]
    \centering
    \begin{tabular}{c|c|c|c|c}
         \textbf{Exp. \#1} & \textbf{Av. Speed} & \textbf{Max Speed} & \textbf{Pr\textsubscript{fail}} & \textbf{Cost}\\
         \hline
         MPPI & 2.46 $\pm$ 0.31 & 3.42 $\pm$ 0.35 & 1.0 &44.1 $\pm$ 7.0\\
         \hline
         Tube-MPPI & 2.37 $\pm$ 0.32 & 2.95 $\pm$ 0.36 & 0.87 & 47.9 $\pm$ 8.2\\
         \hline
         CCSMPPI & 2.33 $\pm$ 0.31 & 3.03 $\pm$ 0.37 & 0.13 & 58.6 $\pm$ 8.5 \\
        \hline
        \hline
         \textbf{Exp. \#2} & \textbf{Av. Speed} & \textbf{Max Speed} & \textbf{Pr\textsubscript{fail}} & \textbf{Cost}\\
         \hline
         MPPI & 1.66 $\pm$ 0.24 & 3.53 $\pm$ 0.32 & 1.0 & 259.6 $\pm$ 73.8 \\
         \hline
         Tube-MPPI & 1.77 $\pm$ 0.26 & 2.85 $\pm$ 0.35 & 0.67 & 95.6 $\pm$ 20.9 \\
         \hline
         CCSMPPI & 1.65 $\pm$ 0.27 & 2.67 $\pm$ 0.29 & 0.07 & 66.1 $\pm$ 8.3
    \end{tabular}
    \caption{Performance Comparision Statistics}
    \label{tab:Performance Comparison}
    \vspace{-0.6cm}
\end{table}

\section{Conclusion}\label{sec:Conclusion}
In this paper, we presented a novel framework for safe trajectory optimization for stochastic linear systems. 
We used Model Predictive Path Integral (MPPI) control for stochastic optimization and a projection-based linearization method for the generation of safe convex regions. 
In addition, we used a Constrained Covariance Steering algorithm 
to safeguard against unmodeled noise disturbances that the MPPI algorithm may not always handle. 
Our numerical simulations have demonstrated that our approach can guarantee safety against unmodeled noise uncertainties as well as unsafe outputs generated by the stochastic optimization algorithm. 
In our future work, we plan to extend our proposed framework to trajectory generation problems for uncertain nonlinear systems based on model-free trajectory optimization algorithms while guaranteeing safety. 

\bibliographystyle{ieeetr}
\bibliography{CSMPPI}

\end{document}